\documentclass[11pt, a4paper]{article}
\usepackage{times}
\usepackage{a4wide}
\usepackage[dvips, hyperindex]{hyperref}
\usepackage[british]{babel}
\usepackage{enumerate}
\usepackage{amsmath, amscd, amsfonts, amssymb, latexsym, theorem, graphicx}
\usepackage[T1]{fontenc}
\usepackage[latin1]{inputenc}

\theoremstyle{change}
{\theorembodyfont{\itshape}   \newtheorem{thm}{Theorem.}[section]}
{\theorembodyfont{\itshape}   \newtheorem{lem}[thm]{Lemma.}}
{\theorembodyfont{\itshape}   \newtheorem{defi}[thm]{Definition.}}
{\theorembodyfont{\itshape}   \newtheorem{rem}[thm]{Remark.}}
{\theorembodyfont{\itshape}   \newtheorem{prop}[thm]{Proposition.}}
{\theorembodyfont{\itshape}   \newtheorem{cor}[thm]{Corollary.}}

\newcommand{\QQ}{\mathbb{Q}}
\newcommand{\CC}{\mathbb{C}}
\newcommand{\ZZ}{\mathbb{Z}}
\newcommand{\FF}{\mathbb{F}}
\newcommand{\NN}{\mathbb{N}}
\newcommand{\Qbar}{{\overline{\QQ}}}

\newcommand{\Fbar}{{\overline{\FF}}}
\newcommand{\rhobar}{{\overline{\rho}}}
\newcommand{\cO}{\mathcal{O}}
\newcommand{\fP}{\mathfrak{P}}
\newcommand{\fQ}{\mathfrak{Q}}
\newcommand{\fp}{\mathfrak{p}}
\newcommand{\GL}{\mathrm{GL}}

\newcommand{\PGL}{\mathrm{PGL}}
\newcommand{\PSL}{\mathrm{PSL}}

\newcommand{\Frob}{\mathrm{Frob}}
\newcommand{\Gal}{\mathrm{Gal}}
\newcommand{\Ind}{\mathrm{Ind}}
\newcommand{\Tr}{\mathrm{tr}}

\newcommand{\triv}{\mathrm{triv}}
\newcommand{\modulo}{\,\mathrm{mod}\,}
\newcommand{\id}{\mathrm{id}}
\newcommand{\proj}{\mathrm{proj}}
\newcommand{\pf}{{\bf Proof. }}
\newcommand{\qed}{\hspace* {.5cm} \hfill $\Box$}
\newcommand{\mat}[4]{
 \left(  \begin{smallmatrix} #1 & #2 \\ #3 & #4 \end{smallmatrix} \right)}
\newcommand{\legendre}[2]{{\left(\frac{#1}{#2}\right)}}

\title{On Modular Forms and the Inverse Galois Problem}
\author{Luis Dieulefait and Gabor Wiese}

\begin{document}

\maketitle

\begin{abstract}
In this article new cases of the Inverse Galois Problem are established.
The main result is that for a fixed integer $n$, there is a positive
density set of primes $p$ such that $\PSL_2(\FF_{p^n})$ occurs as the
Galois group of some finite extension of the rational numbers.
These groups are obtained as projective images of residual modular Galois representations.
Moreover, families of modular forms are constructed such that the images of
all their residual Galois representations are as large as a priori possible.
Both results essentially use Khare's and Wintenberger's notion of good-dihedral primes.
Particular care is taken in order to exclude nontrivial inner twists.

\end{abstract}

\section{Introduction}

The inverse Galois problem asks whether any
given finite group occurs as the Galois group of some Galois extension $K/\QQ$.

For any newform~$f$ and any prime~$\ell$,
the projective image $G$ of the mod~$\ell$ Galois representation
attached to~$f$ is a Galois group over~$\QQ$ by Galois theory.
If $\ell$ is a so-called nonexceptional prime for~$f$, then $G$ is
of the type $\PSL_2(\FF_{\ell^n})$ or $\PGL_2(\FF_{\ell^n})$.
If $\ell$ is exceptional, then $G$ is an abelian group, a dihedral group,
or $A_4$, $S_4$ or~$A_5$.

In this article we construct families of newforms without exceptional primes
(see Theorem~\ref{thm:family}). This is achieved by exploiting the notion of
{\em tamely dihedral primes}, which are essentially the same as the
{\em good-dihedral primes} introduced by Khare and Wintenberger in their proof
of Serre's modularity conjecture~\cite{KW1}.
All the constructed newforms also enjoy the property that they do not have
any nontrivial inner twist and are not CM forms (see Section~\ref{sec:it}).

\begin{figure}\label{fig:results}
\includegraphics[width=8cm]{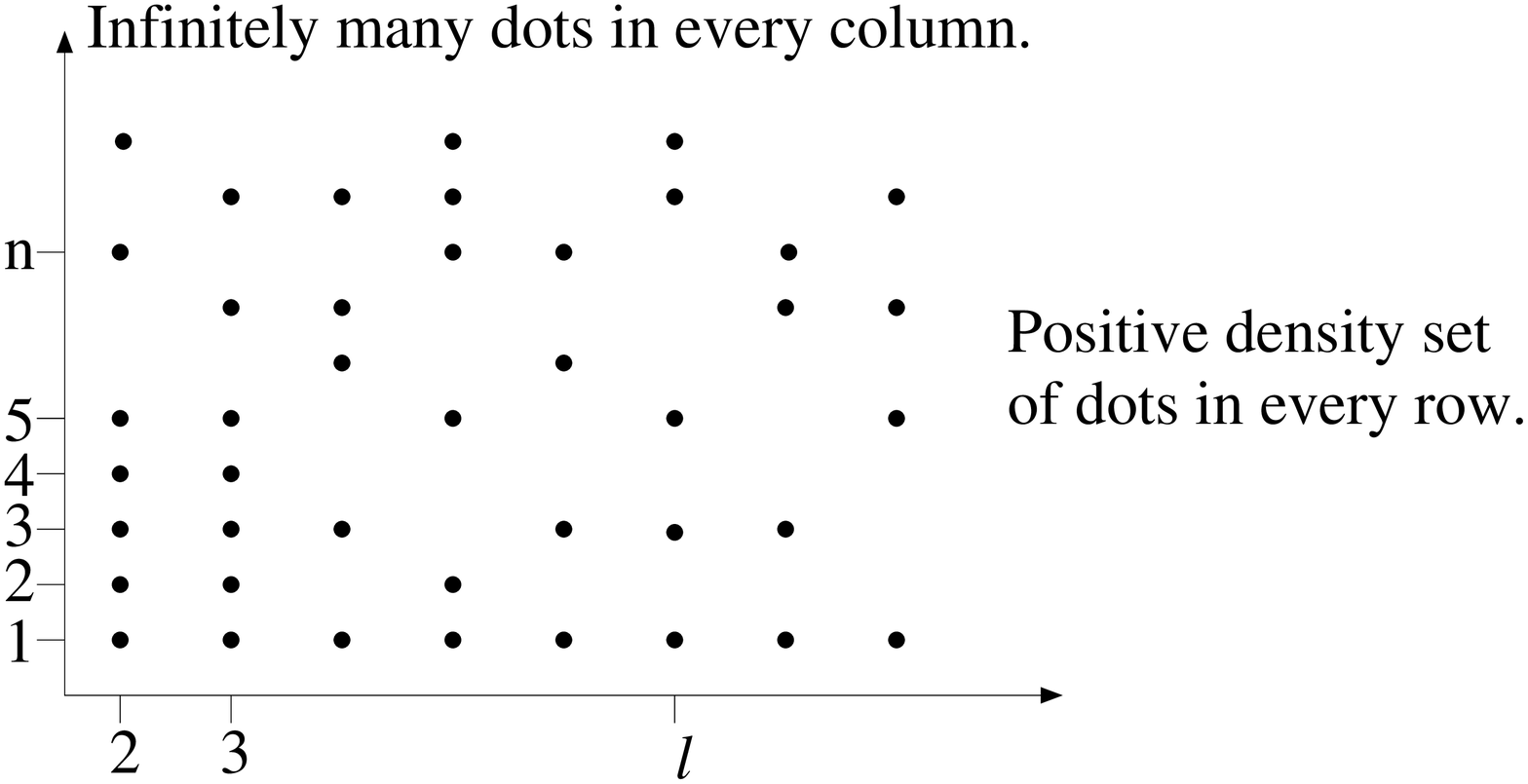}
\caption{Schematic description of the proved cases of $\PSL_2(\FF_{\ell^n})$
occuring as a Galois group over~$\QQ$}
\end{figure}

These techniques, together with some other methods, are applied to the inverse
Galois problem for groups of the type $\PSL_2(\FF_{\ell^n})$ and $\PGL_2(\FF_{\ell^n})$.
Our results on that problem fall in two different categories, which we like to
call the {\em horizontal direction} and the {\em vertical direction}. The terminology
is explained by Figure~\ref{fig:results}.
If a dot exists at position $(\ell,n)$ in the figure, then it is known that
$\PSL_2(\FF_{\ell^n})$ is a Galois group over~$\QQ$.
Theorem~1.1 of~\cite{W} is the best result in the vertical direction to this date.
It says that in every column (i.e.\
for every prime~$\ell$) there are infinitely many dots, i.e.\
there are infinitely many~$n$ such that $\PSL_2(\FF_{\ell^n})$
is a Galois group over~$\QQ$.
Except for finitely many columns, this result was reproved in~\cite{D}
by different methods. We also prove a similar (but slightly weaker)
result in this article (see Remark~\ref{rem:vertical}), but in a uniform way
by using the families without exceptional primes mentioned above.
The vertical result has been vastly generalized to many more groups by
Khare, Larsen and Savin in \cite{KLS} and~\cite{KLS2}.

There were several results for small~$n$ in the horizontal direction
(see \cite{Ri75}, \cite{RV}, \cite{DV}, \cite{Di1}). The main result
of this article is that in every row (i.e.\ for every~$n$) the set
of primes~$\ell$ such that $\PSL_2(\FF_{\ell^n})$ is a Galois group
over~$\QQ$ has a positive density. More precisely, we prove the
following theorem.

\begin{thm}\label{thm:main}
Let $n$ be some integer.
\begin{enumerate}[(a)]
\item There is a positive density set of primes~$\ell$ such
that $\PSL_2(\FF_{\ell^n})$ occurs as the Galois group $\Gal(K/\QQ)$
for a number field~$K$ that ramifies at most at $\ell$, $\infty$ and
two other primes for even~$n$ and three other primes for odd~$n$.
\item Assume that $n$ is odd. There is a positive density set of primes~$\ell$ such
that $\PGL_2(\FF_{\ell^n})$ occurs as the Galois group $\Gal(K/\QQ)$
for a number field~$K$ that ramifies at most at $\ell$, $\infty$ and
two other primes.
\end{enumerate}
\end{thm}

We finish this introduction by pointing out that our method of obtaining
groups of the type $\PSL_2(\FF_{\ell^n})$ or $\PGL_2(\FF_{\ell^n})$
as Galois groups over~$\QQ$ via newforms
is very general: if a group of this type occurs as the Galois group of a
{\em totally imaginary} extension of~$\QQ$, then it is the projective image of the
Galois representation attached to a newform by Serre's modularity conjecture,
which is now a theorem of Khare-Wintenberger (\cite{KW1}, \cite{KW2}, see also \cite{Ki} and
\cite{Di3}).

\begin{prop}
Let $K/\QQ$ be a totally imaginary Galois extension with Galois group $G$ which is either
$\PSL_2(\FF_{\ell^n})$ or $\PGL_2(\FF_{\ell^n})$. Then there exists a modular form
$f$ such that its attached projective mod~$\ell$ Galois representation cuts out the field~$K$.
\end{prop}

\pf
We interpret the number field~$K$ as a projective Galois
representation $\rhobar^\proj: \Gal(\Qbar/\QQ) \to
\PGL_2(\FF_{\ell^n})$. Now we lift the representation to a
continuous representation $\rhobar: \Gal(\Qbar/\QQ) \to
\GL_2(\Fbar_\ell)$ (this is possible, see e.g.~\cite{Q}) and we
observe that the lift is necessarily odd. For, the image of complex
conjugation in $\PGL_2(\FF_{\ell^n})$ being nontrivial, it follows
that its image under $\rhobar$ is a nonscalar involution. As such,
it has determinant $-1$. Finally, one invokes Serre's modularity conjecture
to obtain the modularity.
\qed

\subsection*{Acknowledgements}

Both authors were partially supported by the European Research Training Network
{\em Galois Theory and Explicit Methods} MRTN-CT-2006-035495.
G.~W. also acknowledges partial support by the Sonderforschungsbereich Transregio 45
of the Deutsche Forschungsgemeinschaft.

\subsection*{Notation}

Here we list some notation to be used throughout the article.
We denote by $S_k(N,\chi)$ the $\CC$-vector space of holomorphic
cuspidal modular forms of level~$N$, weight~$k$ and Dirichlet charater~$\chi$.
By $\chi_\triv$ we mean the trivial Dirichlet character; instead
of $S_k(N,\chi_\triv)$ we also write $S_k(\Gamma_0(N))$.
If $K$ is a number field, we denote by $\cO_K$ its ring of integers,
and by $\cO_{K,\Lambda}$ its completion at a maximal ideal $\Lambda \lhd \cO_K$.
For a prime~$q$, we let $D_q$ stand for $\Gal(\Qbar_q/\QQ_q)$ and denote
by $I_q \subset D_q$ the inertia group. By $W_q$ we denote the Weil group of~$\QQ_q$.
By $\zeta_n$ for an integer~$n$ we always denote a primitive $n$-th root of unity.

More notation will be introduced in the text.

\section{Complex Multiplication and inner twists}\label{sec:it}

In this section we review essential facts on complex multiplication and inner twists.
Let $f$ be some cuspidal modular form
of level~$N$, weight~$k$ and Dirichlet character~$\chi$ with $q$-expansion
$\sum_{n \ge 1} a_n(f) q^n$ (as usual, $q = q(\tau) = e^{2 \pi i \tau}$).
The coefficient field of~$f$ is defined as $\QQ_f = \QQ(a_n(f) \,:\,(n,N)=1)$.
It has the natural subfield $F_f = \QQ\left(\frac{a_n(f)^2}{\chi(n)} \,:\,(n,N)=1\right)$,
which we call the {\em twist invariant coefficient field of~$f$}, since
it is invariant under replacing the modular form~$f$ by any of its twists.
The behaviour of the Hecke operators under the Petersson scalar product
yields the formula
\begin{equation}\label{eqcc}
   \overline{a_p(f)} = \chi(p)^{-1} a_p(f),
\end{equation}
whence $\frac{a_p(f)^2}{\chi(p)} = |a_p(f)|^2$. Thus, $F_f$ is totally real.
It is well known that $\QQ_f$ is either a CM field or totally real.
In particular, if $f$ has trivial nebentype the latter case occurs by
Equation~\ref{eqcc}.

The modular form $f$ is said to have {\em complex multiplication (CM)} if there exists
a Dirichlet character~$\epsilon$ such that
\begin{equation}\label{eqcm}
  a_p (f \otimes \epsilon) = a_p(f) \epsilon(p) = a_p(f)
\end{equation}
for almost all primes~$p$ (i.e.\ all but finitely many).

A twist of~$f$ by a Dirichlet character~$\epsilon$ is said to be {\em inner}
if there exists a field automorphism $\sigma_\epsilon: \QQ_f \to \QQ_f$
such that
\begin{equation}\label{eqin}
  a_p (f \otimes \epsilon) = a_p(f) \epsilon(p) = \sigma_\epsilon (a_p(f))
\end{equation}
for almost all primes~$p$.

For a discussion of inner twists we refer the reader to \cite{Rib80} and
\cite{Rib85}. Here we content ourselves by listing some statements that will be needed for
the sequel.
For the rest of this section we suppose
that $f$ does not have CM.
Under this assumption $\sigma_\epsilon$ is unique (if it exists).
Moreover, for two distinct inner twists $\epsilon$ and $\delta$, the
field automorphism $\sigma_\epsilon$ is different from $\sigma_\delta$.

The inner twists and the Dirichlet character of~$f$ satisfy the relation
\begin{equation}\label{eqcharacter}
   \chi \epsilon^2 = \sigma_{\epsilon}(\chi).
\end{equation}
If $\chi$ takes only real values (i.e.\ is trivial or quadratic), it follows
that $\epsilon^2=1$.
This in turn implies
$$ \sigma_\epsilon^2(a_p(f)) = \sigma_\epsilon (a_p(f) \epsilon(p))
= \sigma_\epsilon(a_p(f)) \epsilon(p) = a_p(f) \epsilon(p)^2 = a_p(f),$$
whence $\sigma_\epsilon^2 = 1$.

The $\sigma_\epsilon$ belonging to the inner twists
of~$f$ form an abelian subgroup~$\Gamma$ of the automorphism group of~$\QQ_f$.
The field $F_f$ is the subfield of~$\QQ_f$ fixed by~$\Gamma$.
If the nebentype of~$f$ only takes real values,
it follows that $\Gamma$ is an elementary abelian $2$-group and
thus that $[\QQ_f:F_f]$ is a power of~$2$.

We define a number field $K_\Gamma$ as follows. Consider the inner twists
$\epsilon_1,\dots,\epsilon_r$ as
characters of the absolute Galois group $\Gal(\Qbar/\QQ)$.
Let $K_\Gamma$ be the minimal number field on which all $\epsilon_i$
for $1\le i \le r$ are trivial, i.e.\ the field such that its absolute
Galois group is the kernel of the map
$\Gal(\Qbar/\QQ) \xrightarrow{\epsilon_1,\dots,\epsilon_r} \CC^\times \times \dots \times \CC^\times$.

\section{On the images of modular Galois representations}

Let, as before, $f = \sum_{n=1}^\infty a_n q^n$ (with $a_n \in \CC$)
be a cuspidal modular form of level~$N$, weight~$k$ and
Dirichlet character~$\chi$. We now assume that $f$ is a Hecke eigenform.
Also as before, let $\QQ_f$ be the coefficient field of~$f$; it is naturally
a subfield of~$\CC$.
By a construction of Shimura and Deligne and the local Langlands
correspondence for $\GL_2$ one can
attach to~$f$ a $2$-dimensional compatible system of Galois representations
$(\rho_{f,\iota})$ of $\Gal(\Qbar/\QQ)$. We now describe this system, following
Khare and Wintenberger~\cite{KW0}.

The compatible system $(\rho_{f,\iota})$ consists of the data of:
\begin{enumerate}[(i)]
\item for each rational prime~$\ell$ and each embedding $\iota: \QQ_f \hookrightarrow \Qbar_\ell$
a continuous semisimple representation $\rho_{f,\iota}: \Gal(\Qbar/\QQ) \to \GL_2(\Qbar_\ell)$,
\item for each rational prime $q$ a Frobenius semisimple Weil-Deligne representation $r_q$
with values in $\GL_2(\QQ_f)$ such that
\begin{enumerate}[(a)]
\item $r_q$ is unramified for all~$q$ outside a finite set,
\item for each rational prime~$\ell$, for each rational prime~$q \neq \ell$ and
for each $\iota: \QQ_f \hookrightarrow \Qbar_\ell$, the Frobenius semisimple Weil-Deligne
representation associated to $\rho_{f,\iota}|D_q$ is conjugate
to~$r_q$ via the embedding~$\iota$.
\end{enumerate}
\item a third condition concerning $q=\ell$ (which we do not need).
\end{enumerate}
The system $\rho_{f,\iota}$ is attached to~$f$ in the sense that
for each rational prime $q \nmid N\ell$, the characteristic polynomial of
$\rho_{f,\iota}(\Frob_q)$ is $X^2 - \iota(a_q)X + \iota(\chi(q) q^{k-1})$.

We also introduce a different description of the $\rho_{f,\iota}$, which we will often
use below. Let $\iota: \QQ_f \to \Qbar_\ell$ be an embedding. Denote by
$\QQ_{f,\iota}$ the closure of $\iota(\QQ_f)$ and by
$\cO_{f,\iota}$ the closure of $\iota(\cO_{\QQ_f})$ in~$\Qbar_\ell$. Let $(\pi)$ be
the maximal ideal of the local complete discrete valuation ring~$\cO_{f,\iota}$.
Then $\Lambda := \cO_{\QQ_f} \cap \iota^{-1}\big((\pi)\big)$ is a maximal ideal
of~$\cO_{\QQ_f}$ above~$\ell$ and $\QQ_{f,\iota}$ can be identified with
$\QQ_{f,\Lambda}$, the completion of $\QQ_f$ at~$\Lambda$. To simplify notation,
we denote by $\cO_{f,\Lambda}$ the integers of~$\QQ_{f,\Lambda}$.
Thus, we can identify $\rho_{f,\iota}$ with
$$ \rho_{f,\Lambda}: \Gal(\Qbar/\QQ) \to \GL_2(\cO_{f,\Lambda}).$$
More precisely, the composition of $\rho_{f,\Lambda}$ with the natural inclusion
$\cO_{f,\Lambda} \hookrightarrow \Qbar_\ell$ equals~$\rho_{f,\iota}$.

In this article, we shall mostly be concerned with the reduction of these
representations 'modulo~$\ell$', by which we mean the following.
Let $\rhobar_{f,\Lambda}$ be the semisimplification of the
reduction of~$\rho_{f,\Lambda}$ modulo~$\Lambda$
and $\rhobar_{f,\Lambda}^\proj$ its projective quotient, i.e.\
$\rhobar_{f,\Lambda}$ composed with the natural projection
$\GL_2(\FF_\Lambda) \twoheadrightarrow \PGL_2(\FF_\Lambda)$, where we write
$\FF_\Lambda = \cO_{\QQ_f}/\Lambda$ for the residue field of~$\Lambda$.

\begin{thm}[Ribet]\label{thm:ribet}
Suppose that $f$ does not have CM and that $k \ge 2$.
Then for almost all maximal ideals $\Lambda$ of $\cO_{\QQ_f}$
(equivalently, for almost all~$\iota$ as above), the image
$\rhobar_{f,\Lambda}(\Gal(\Qbar/K_\Gamma))$ is conjugate to
  $$ \{g \in \GL_2(\FF_\lambda) \,:\, \det(g) \in \FF_\ell^{\times (k-1)} \},$$
where $K_\Gamma$ is the field defined in Section~\ref{sec:it},
$\FF_\ell$ is the prime field of~$\FF_\Lambda$
and $\FF_\lambda$ is the residue field of $F_f$ at $\lambda = \Lambda \cap F_f$,
i.e.\  $\FF_\lambda = \cO_{F_f}/\lambda$.
\end{thm}

\pf
It suffices to take Ribet~\cite[Thm.~3.1]{Rib85} mod~$\Lambda$.
\qed

\begin{cor}\label{cor:ribet}
Under the assumptions of Theorem~\ref{thm:ribet}, for almost all~$\Lambda$ the image
$\rhobar_{f,\Lambda}^\proj(\Gal(\Qbar/\QQ))$ is either $\PSL_2(\FF_\lambda)$
or $\PGL_2(\FF_\lambda)$.
\end{cor}

\pf
Let $H := \rhobar_{f,\Lambda}(\Gal(\Qbar/K_\Gamma))$. It is a normal subgroup
of $G := \rhobar_{f,\Lambda}(\Gal(\Qbar/\QQ))$. Consequently,
$H/(H \cap \Fbar_\ell^\times)$ is a normal subgroup of $G/(G \cap \Fbar_\ell^\times)$.
By the classification of finite subgroups of $\PGL_2(\Fbar_\ell)$, it follows
for almost all~$\Lambda$ that $H$ is either $\PSL_2(\FF_\lambda)$ or $\PGL_2(\FF_\lambda)$ and
that $G$ is either $\PSL_2(\FF)$ or $\PGL_2(\FF)$ for some extension $\FF/\FF_\lambda$.
As $\PSL_2(\FF)$ is simple (for $\#\FF \ge 5$), $\FF = \FF_\lambda$ follows.
\qed

\begin{defi}
Keep the assumptions of Theorem~\ref{thm:ribet}. A maximal ideal $\Lambda$ of $\cO_{\QQ_f}$
is called {\em exceptional} if $\rhobar_{f,\Lambda}^\proj(\Gal(\Qbar/\QQ))$ is
neither $\PSL_2(\FF_\lambda)$ nor $\PGL_2(\FF_\lambda)$.
\end{defi}

By Corollary~\ref{cor:ribet}, every form~$f$ without CM only has finitely many exceptional
primes. The classification of finite subgroups of $\PGL_2(\Fbar_\ell)$ yields that
$\rhobar_{f,\Lambda}^\proj(\Gal(\Qbar/\QQ))$ is either an abelian group, a dihedral group,
$A_4$, $S_4$ or $A_5$ for an exceptional prime~$\Lambda$.

\section{Tamely dihedral representations}

It was observed by Khare and Wintenberger (\cite{KW1}) that exceptional primes
smaller than a given bound can be avoided
by imposing a certain local ramification behaviour of the Galois representation.
We shall use this idea in order to exclude CM and inner twists.

We will formulate the crucial definition in terms of Weil-Deligne representations.
As a reference for these, we use \cite{T}, Section~4.
Let $K$ be a finite extension of~$\QQ_q$ (with $q$ a prime)
and $E$ a number field (later $E=\QQ_f$).
A $2$-dimensional Weil-Deligne representation of~$K$
with values in~$E$ can be described as a pair $(\widetilde{\rho},\widetilde{N})$
(we put tildes in order to avoid any possible confusion with Galois representations
and levels of modular forms), where
$\widetilde{\rho}: W_K \to \GL_2(E)$ is a continuous representation of the Weil group of~$K$
for the discrete topology on~$\GL_2(E)$ and $\widetilde{N}$ is a nilpotent endomorphism
of~$E^2$ satisfying a certain commutativity relation with~$\widetilde{\rho}$.

\begin{defi}\label{defi:td}
Let $\QQ_{q^2}$ be the unique unramified degree~$2$ extension of~$\QQ_q$.
Denote by $W_q$ and $W_{q^2}$ the Weil group of~$\QQ_q$ and $\QQ_{q^2}$, respectively.

A $2$-dimensional Weil-Deligne representation $r_q=(\widetilde{\rho},\widetilde{N})$ of $\QQ_q$
with values in~$E$ is called {\em tamely dihedral of order~$n$} if $\widetilde{N}=0$ and
there is a tame character $\psi: W_{q^2} \to E^\times$
whose restriction to the inertia group $I_q$ (which is naturally
a subgroup of~$W_{q^2}$) is of niveau~$2$
(i.e.\ it factors over $\FF_{q^2}^\times$ and not over $\FF_q^\times$) and of order~$n>2$
such that $\widetilde{\rho} \cong \Ind_{W_{q^2}}^{W_q} (\psi)$.

We say that a Hecke eigenform~$f$ is {\em tamely dihedral of order~$n$} at the
prime~$q$ if the Weil-Deligne representation $r_q$ at~$q$
belonging to the compatible system $(\rho_{f,\iota})$ is tamely dihedral of order~$n$.
\end{defi}

If the compatible system $(\rho_{f,\iota})$ is tamely dihedral of order~$n$ at~$q$,
then (e.g.\ by \cite{T}, 4.2.1)
for all $\iota: \QQ_f \to \Qbar_\ell$ with $\ell \neq q$, the restriction of
$\rho_{f,\iota}$ to $D_q$ is of the form
$\Ind_{\Gal(\Qbar_q/\QQ_{q^2})}^{\Gal(\Qbar/\QQ_q)} (\iota \circ \psi)$.
The point in our application is that
the reduction modulo~$\ell$ is of the very same form,
i.e.\ if $\overline{\psi}_\Lambda$ denotes the reduction of~$\psi$
modulo~$\Lambda$, which is a character of the same order if $\Lambda$ and
$n$ are coprime, then
$$\rhobar_{f,\Lambda}|_{D_q} =
\Ind_{\Gal(\Qbar_q/\QQ_{q^2})}^{\Gal(\Qbar/\QQ_q)} (\overline{\psi}_\Lambda).$$
Moreover, if $n = p^r$ for some odd prime~$p$, then $q \equiv -1 \modulo p$, since
the character is of niveau~$2$.
Conversely, if we take $\psi$ to be a totally ramified character of $\Gal(\Qbar_q/\QQ_{q^2})$
of order~$n$ such that $n$ divides $q+1$ and $(n,q(q-1))=1$, then this automatically
ensures that $\psi$ is of niveau~$2$.

We also mention that the image of $\widetilde{\rho}$ as in Definition~\ref{defi:td}
is isomorphic as an abstract group to the dihedral group~$D_n$ of order~$2n$
if $\psi$ is totally ramified and $(n, q(q-1))=1$.
This is due to the fact that this condition forces the
restriction of the determinant to~$\Gal(\Qbar_q/\QQ_{q^2})$ to be trivial.

A tamely dihedral prime~$q$ is called a {\em good-dihedral prime}
by Khare and Wintenberger \cite{KW1} if some additional properties
(of a different nature) are satisfied. Good-dihedral primes play an
important role in Khare and Wintenberger's proof of Serre's conjecture.

We now collect some very simple lemmas that will afterwards be applied in order
to exclude nontrivial inner twists.

\begin{lem}\label{lem:td}
Let $K$ be a topological field, $q$ a prime and $n>2$ an integer coprime to~$q(q-1)$.
Let $\epsilon: \Gal(\Qbar_q/\QQ_q) \to K^\times$
and $\psi, \psi': \Gal(\Qbar/\QQ_{q^2}) \to K^\times$ be characters. Suppose that
$\psi$ and $\psi'$ are both of order~$n$.

If $\Ind_{\QQ_{q^2}}^{\QQ_q} (\psi) \cong \Ind_{\QQ_{q^2}}^{\QQ_q} (\psi') \otimes \epsilon$,
then $\epsilon$ is unramified.
\end{lem}

\pf
It follows directly that the order of $\epsilon|_{\Gal(\Qbar_q/\QQ_{q^2})}$ divides~$n$.
If $\epsilon$ were ramified, then the order of $\epsilon$ restricted to the inertia group at~$q$
would divide $q-1$ times a power of~$q$. But, this is precisely excluded in the assumption.
\qed

\begin{lem}\label{lem:oddorder}
Let $K$ be a topological field.
Let $q>2$ be a prime, $\epsilon: \Gal(\Qbar_q/\QQ_q) \to K^\times$ a quadratic character
and $\rho, \rho': \Gal(\Qbar_q/\QQ_q) \to \GL_2(K)$ representations such that
$\rho' \cong \rho \otimes \epsilon$.
If $\rho(I_q)$ and $\rho'(I_q)$ can be conjugated to lie in the upper triangular
matrices such that all elements on the diagonal have odd order,
then $\epsilon$ is unramified.
\end{lem}

\pf
The oddness of the order of the elements on the diagonal forces
the quadratic character~$\epsilon$ to be unramified.
\qed

\begin{lem}\label{lem:sqf}
Let $K$ be a topological field.
Let $q$ be a prime, $\epsilon: \Gal(\Qbar_q/\QQ_q) \to K^\times$ a character
and $\rho: \Gal(\Qbar_q/\QQ_q) \to \GL_2(K)$ be a representation.
If the conductors of $\rho$ and of $\rho \otimes \epsilon$ both divide~$q$, then
$\epsilon$ or $ \epsilon \det(\rho)$ is unramified.
\end{lem}

\pf
By the definition of the conductor, $\rho$ restricted to the inertia group $I_q$ is of the
form $\mat 1*0 \delta$ with $\delta = \det(\rho)|_{I_q}$. Consequently,
the restriction to~$I_q$ of $\rho \otimes \epsilon$ looks like $\mat \epsilon*0 {\epsilon \delta}$.
Again by the definition of the conductor, either the restriction of $\epsilon$ to $I_q$ is trivial
or the restriction of $\epsilon \delta$ is.
\qed
\medskip

Our main result for controlling inner twists and CM is the following theorem.

\begin{thm}\label{thm:it}
Let $f \in S_k(N,\chi)$ be a normalized Hecke eigenform.
\begin{enumerate}[(a)]
\item If $\chi$ only takes real values,
then any inner twist of~$f$ is at most tamely ramified at any odd prime.
\item Suppose $\chi=\chi_\triv$.
Let $q \mid N$ be a prime such that for some $\iota: \QQ_f \to \Qbar_\ell$ with $\ell \neq q$
the image $\rho_{f,\iota}(I_q)$ can be conjugated to lie in the upper triangular
matrices such that the diagonal elements have odd order.
Then any inner twist of~$f$ is unramified at~$q$.
\item Let $q \mid\mid N$ be a prime and suppose that $\chi$ is unramified at~$q$. Then
any inner twist of~$f$ is unramified at~$q$.
\item Let $q$ be a prime such that $q^2 \mid\mid N$ and $f$ is tamely dihedral at~$q$
of odd order~$n \ge 3$ such that $n$ and $q(q-1)$ are coprime.
Then any inner twist of~$f$ is unramified at~$q$.
\end{enumerate}
\end{thm}

\pf
(a) If $\chi$ only takes real values, any inner twist is necessarily quadratic
(see Section~\ref{sec:it}), whence it is at most tamely ramified away from~$2$.

For (a-c) we let $\epsilon$ be an inner twist with corresponding field automorphism
$\sigma: \QQ_f \to \QQ_f$.
We include the case $\sigma=\id$, when $\epsilon$ comes from CM.
For all $\iota: \QQ_f \to \Qbar_\ell$, we then have
$$\rho_{f,\iota} \otimes \epsilon \cong \rho_{f,\iota \circ \sigma},$$
since the traces of any Frobenius element at any unramified prime~$p$ are equal:
$$\Tr((\rho_{f,\iota} \otimes \epsilon)(\Frob_p)) = \iota (a_p(f) \epsilon(p))
= \iota(\sigma(a_p(f))) = \Tr(\rho_{f,\iota\circ \sigma}(\Frob_p)).$$

(b) If $\rho_{f,\iota}(I_q)$ can be conjugated into the upper triangular matrices
such that the diagonal elements have odd order, then it follows that
the Weil-Deligne representation $r_q = (\widetilde{\rho},\widetilde{N})$
(in the terminology of \cite{T},~4.1.2) is such that $\widetilde{\rho}$ also can be
conjugated into the upper triangular matrices with only odd order elements
on the diagonal. Hence, this property also holds for
$\rho_{f,\iota \circ \sigma}(I_q)$. Consequently,
Lemma~\ref{lem:oddorder} yields the statement, as $\epsilon$ again has to be
at most quadratic in this case.

(c) The conductors at~$q$ of both $\rho_{f,\iota}$ and $\rho_{f,\iota \circ \sigma}$ divide~$q$.
From Lemma~\ref{lem:sqf} it hence follows that $\epsilon$ is unramified at~$q$,
since the determinant of the representation is unramified at~$q$.

(d) If $r_q$ is tamely dihedral of order~$n$ as in the assumption,
the restriction of $\rho_{f,\iota}$ to $D_q$ is of the form
$\Ind_{\Gal(\Qbar/\QQ_{q^2})}^{\Gal(\Qbar/\QQ_q)} (\iota \circ \psi)$,
and similarly for $\rho_{f,\iota\circ\sigma}$.
By Lemma~\ref{lem:td}, $\epsilon$ is unramified at~$q$.
\qed

\begin{cor}\label{cor:noit}
Let $f \in S_k(N,\chi_\triv)$ be a Hecke eigenform such that for all primes $q \mid N$
\begin{itemize}
\item $q \mid\mid N$, or
\item $q^2 \mid\mid N$ and $f$ is tamely dihedral at~$q$ of some order $n \ge 3$
such that $(n,q(q-1))=1$, or
\item $\rho_{f,\iota}(I_q)$ can be conjugated into the upper triangular matrices
such that all the elements on the diagonal have odd order for some $\iota: \QQ_f \to \Qbar_\ell$
with $\ell \neq q$.
\end{itemize}

Then $f$ does not have any nontrivial inner twists and no CM.
\end{cor}

\pf
By Theorem~\ref{thm:it}, any inner twist (or character corresponding to CM)
is everywhere unramified, hence trivial.
\qed
\medskip

The next proposition will be essential in the application of modular forms to
the inverse Galois problem in the horizontal direction.

\begin{prop}\label{prop:ordp}
Let $f \in S_k(Nq^2,\chi)$ be a Hecke eigenform which is tamely dihedral of
some order~$n>2$ at~$q$ such that $(n,q(q-1))=1$.
Then $F_f$ contains $\QQ(\zeta_{n}+\overline{\zeta_{n}})$.
\end{prop}

\pf
Let $K_\Gamma$ be the field defined in Section~\ref{sec:it}.
From Theorem~\ref{thm:it} it follows that all inner twists are unramified at~$q$.
Hence, the inertia group $I_q$ can be considered as a subgroup of
$\Gal(\Qbar/K_\Gamma)$.
Let $\Lambda$ be any prime of $\cO_{\QQ_f}$ not dividing~$nq$.
By assumption, the image $\rho_{f,\Lambda}(I_q)$ contains an element
of the form $\mat {\zeta_n} * 0 {\overline{\zeta_n}}$.
By Theorem~\ref{thm:ribet}, $F_{f,\lambda}$ contains its
trace, i.e.\ $\zeta_n + \overline{\zeta_n}$, where $\lambda = \Lambda \cap F_f$.
This immediately implies that the field extension
$F_f(\zeta_n + \overline{\zeta_n})/F_f$
is of degree~$1$, as almost all primes are completely split in it.

Alternatively, one could also derive this proposition by similar arguments
directly from the Weil-Deligne representation $r_q$ at~$q$.
\qed

\section{Constructions of newforms that are tamely dihedral at some prime}

In this section we collect tools for constructing newforms
and consequently Galois representations which are tamely dihedral at some prime.

\begin{thm}[Diamond, Taylor: {\em Level Raising}]\label{thm:dt}
Let $N \in \NN$, $k \ge 2$ and let $p > k+1$ be a prime not dividing~$N$.
Let $f \in S_k(N,\chi)$ be a newform such that $\rhobar_{f,\fP}$
is irreducible with $\cO_{\QQ_f} \rhd \fP \mid p$.
Let, furthermore, $q \nmid N$ be a prime such that $q \equiv -1 \modulo p$
and $\Tr(\rhobar_{f,\fP}(\Frob_q))=0$.

Then there exists a newform $g \in S_k(Nq^2,\tilde{\chi})$ such that
\begin{enumerate}[(i)]
\item $\rhobar_{g,\fp} \cong \rhobar_{f,\fP}$ for some prime $\fp \mid p$ of $\cO_{\QQ_g}$.

\item \label{part:gf} $g$ is tamely dihedral of order~$p^r$ for some $r > 0$ at~$q$.
\end{enumerate}
\end{thm}

\pf
This is easily deduced from Theorem~A of~\cite{DT}, using the local
Langlands correspondence. For details see Corollary 2.6 of~\cite{W}.
\qed
\medskip

A very elaborate possibility of prescribing a tamely dihedral prime
is Theorem~5.1.4 of~\cite{KW1}, which does not need the modularity assumption.

To our knowledge, the best result for prescribing the local behaviour at ramified primes
for modular Galois representations is the following theorem by Jared Weinstein, which
actually holds for Hilbert modular forms. We only formulate it over~$\QQ$.
A global inertial type $\tau$ is a collection of local inertial types $(\tau_\nu)_\nu$
(or equivalently - by the local Langlands correspondence - inertial Weil-Deligne parameters)
for $\nu$ running through the places of~$\QQ$. For the precise
definitions we refer to~\cite{We}. We just say that the local inertial type at
a finite place~$\nu$ determines the restriction of any $\Lambda$-adic representation
to the inertia group at~$\nu$ uniquely for $\cO_{\QQ_f} \rhd \Lambda \nmid \nu$.
Any newform $f$ uniquely determines a global inertial type, which we denote by $\tau(f)$.

\begin{thm}[Weinstein]\label{thm:weinstein}
Up to twisting by one-dimensional characters, the set of global inertial types $\tau$
for which there is no newform~$f$ with $\tau = \tau(f)$ is finite.
\end{thm}

\pf
\cite{We}, Corollary~1.2.
\qed
\medskip

This means that by making the weight big enough, for any global inertial type~$\tau$ there is
some newform~$f$ with $\tau = \tau(f)$. Alternatively, there is always some newform~$f$
of a chosen weight with $\tau=\tau(f)$ if enough primes ramify.
Weinstein's result is extremely strong and could be used in several places in this
article. We, however, chose more classical arguments like level raising.

For the construction of tamely dihedral modular forms via level raising we need
the following lemma.

\begin{lem}\label{lem:raise}
Let $p_1,\dots,p_r$ be primes and let $p$ be a prime different
from all the $p_i$ such that $p \equiv 1 \modulo 4$.
Let $\rhobar_p^\proj: \Gal(\Qbar/\QQ)\to \PGL_2(\Fbar_p)$ be an odd Galois representation
with image $\PSL_2(\FF_{p^s})$ or $\PGL_2(\FF_{p^s})$ such that the image of any
complex conjugation is contained in $\PSL_2(\FF_{p^s})$.

The set of primes~$q$ such that
\begin{enumerate}[(i)]
\itemsep=0cm plus 0pt minus 0pt
\item $q \equiv p-1 \modulo p^2$,
\item $q$ splits in $\QQ(i,\sqrt{p_1},\dots,\sqrt{p_r})$ and
\item $\rhobar_p^\proj(\Frob_q)$ lies in the same conjugacy class
as $\rhobar_p^\proj(c)$, where $c$ is any complex conjugation,
\end{enumerate}
has a positive density.
\end{lem}

\pf
The proof is adapted from \cite{KW1}, Lemma~8.2, and closely follows the
proof of Lemma~3.2 of~\cite{W}.
Let $L := \QQ(\zeta_{p^2},i,\sqrt{p_1},\dots,\sqrt{p_r})$
and let $K/\QQ$ be such that
$\Gal(\Qbar/K) = \ker(\rhobar_p^\proj)$. Conditions (i) and~(ii) must
be imposed on the field~$L$ and Condition~(iii) on~$K$.
If $\Gal(K/\QQ)$ is $\PSL_2(\FF_{p^s})$, then
the lemma follows directly from Chebotarev's density theorem,
as the intersection $L \cap K$ is~$\QQ$, since $\PSL_2(\FF_{p^s})$
is a simple group. If $\Gal(K/\QQ)$ is $\PGL_2(\FF_{p^s})$, the intersection $L \cap K = M$
is either trivial or an extension of~$\QQ$ of degree~$2$.
By assumption the image of any complex conjugation lies in $\Gal(K/M) \cong \PSL_2(\FF_{p^s})$.
Hence any~$q$ satisfying Condition~(iii) is split
in $M/\QQ$. As $p \equiv 1 \modulo 4$, complex conjugation fixes
the quadratic subfield of~$\QQ(\zeta_{p^2})$, whence any prime~$q$ satisfying
Conditions (i) and~(ii) is also split in $M/\QQ$. Hence,
we may again appeal to Chebotarev's density theorem, proving
the lemma.
\qed
\medskip

In the next proposition we show that in certain cases we can 'add' tamely dihedral primes
to newforms in such a way that all the local behaviours at the primes dividing
the conductor remain essentially the same.
The idea behind this proposition is that for a given newform~$f$ we choose a newform~$g$
such that the compatible systems of Galois representations of $f$ and~$g$ are linked
mod~$p$ for a prime~$p$ that is big enough to preserve all essential local properties.

\begin{prop}\label{prop:tamedihedral}
Let $f \in S_k(N,\chi_\triv)$ be a newform with odd~$N$ such that for all $\ell \mid N$
\begin{enumerate}[(i)]
\item\label{eins} $\ell \mid\mid N$ or
\item\label{zwei} $\ell^2 \mid\mid N$ and $f$ is tamely dihedral at~$\ell$ of order~$n_\ell>2$ or
\item\label{drei} $\ell^2 \mid N$ and $\rho_{f,\iota}(D_\ell)$ can be conjugated to lie in the
upper triangular matrices such that the elements on the diagonal all have odd order
for some $\iota: \QQ_f \hookrightarrow \Qbar_s$ with a primes $s \neq \ell$.
\end{enumerate}
Let $\{p_1,\dots,p_r\}$ be any finite set of primes.

Then for almost all primes $p \equiv 1 \modulo 4$ there is a set~$S$ of primes
of positive density which are completely split in $\QQ(i,\sqrt{p_1},\dots,\sqrt{p_r})$
such that for all $q \in S$ there is a newform
$g \in S_k(Nq^2,\chi_\triv)$ which is tamely dihedral at~$q$ of order~$p$ and
for all~$\ell \mid N$ we have
\begin{enumerate}
\item[(\ref{zwei})] $\ell^2 \mid\mid N$ and $g$ is tamely dihedral at~$\ell$
of order~$n_\ell>2$ or
\item[(\ref{drei})] $\rho_{g,\iota}(D_\ell)$ can be conjugated to lie in the
upper triangular matrices such that the elements on the diagonal all have odd order
for some $\iota: \QQ_g \hookrightarrow \Qbar_s$ with a prime $s \neq \ell$.
\end{enumerate}
Moreover, $f$ and $g$ do not have any nontrivial inner twists and no CM.
\end{prop}

\pf
By Theorem~\ref{thm:it}, $f$ does not have any nontrivial inner twists and no CM.
For $p$ we may choose any prime $p \equiv 1 \modulo 4$ which is bigger than~$N$,
bigger than $k+1$, coprime to all~$n_\ell$ and such that there is $\fP \mid p$ with the property
that $\rhobar_{f,\fP}$ is irreducible (see Corollary~\ref{cor:ribet}).

As $S$ we take the set provided by Lemma~\ref{lem:raise}. Note that the assumptions
of the lemma are satisfied, as complex conjugation necessarily lies in~$\PSL_2$,
as there are no nontrivial inner twists and $-1$ is a square in $\FF_p^\times$.

For any $q \in S$, Theorem~\ref{thm:dt} provides us with a newform
$g \in S_k(Nq^2,\chi)$ which is tamely dihedral at~$q$ of order~$p^r>1$.
In fact, $r=1$ and $\chi=\chi_\triv$. For, as $p^2$ does not divide $q^2-1$,
it follows that there is no niveau~$1$ or niveau~$2$
character of the inertia group~$I_q$ of order~$p^2$.
That $\chi$ is unramified at~$q$ is clear, since
the determinant of the restriction
to inertia at~$q$ is $\psi \psi^q = \psi \psi^{-1} = 1$ (see also the
discussion following Definition~\ref{defi:td}).
Hence, $\chi$ is a character $(\ZZ/N\ZZ)^\times \to \CC^\times$.
As there is a prime $\cO_{\QQ_g} \rhd \fp \mid p$ such that
$\rhobar_{f,\fP} \cong \rhobar_{g,\fp}$, the order of~$\chi$ is a power of~$p$.
As $p$ is bigger than~$N$, it can only be the trivial character.

Next we check that conditions (\ref{zwei}) and (\ref{drei}) persist for~$g$.
This follows again from $\rhobar_{f,\fP} \cong \rhobar_{g,\fp}$, since $p$
is 'big enough'. More precisely, we argue as follows.
Let $\ell \mid N$ be a prime and denote by $(\widetilde{\rho}_f,\widetilde{N}_f)$
be the Weil-Deligne representation attached to~$f$ at~$\ell$.

We now assume that the order of $\widetilde{\rho}_f(I_\ell)$ is divisible by
an odd prime. Note that this condition is satisfied at all primes~$\ell$ in
(\ref{zwei}) and~(\ref{drei}).
We first claim that $\rho_{f,\fP}|_{D_\ell}$ is irreducible if and only if
$\rhobar_{f,\fP}|_{D_\ell}$ is irreducible.
Since the other direction is trivial, we now assume that $\rho_{f,\fP}|_{D_\ell}$
is irreducible. This implies that the representation $\widetilde{\rho}_f$ of the
Weil group of~$\QQ_\ell$ is also irreducible and consequently $\widetilde{N}_f = 0$.
It follows that the projectivization of $\rho_{f,\fP}|_{D_\ell}$ is
$\Ind_{K}^{\QQ_\ell}(\psi)$ with $K/\QQ_\ell$ of degree~$2$ and $\psi$ a nontrivial
character of $\Gal(\Qbar_\ell/K)$ which is different from its conjugate by the
nontrivial element in $\Gal(K/\QQ_\ell)$. Our assumption now means that
the order of~$\psi$ restricted to~$I_\ell$ is divisible by an odd prime. Local
class field theory moreover yields that this prime divides $\ell(\ell^2-1)$ and
is thus different from~$p$, since $p > N$. This implies that the
projectivization of $\rhobar_{f,\fP}|_{D_\ell}$,
which is $\Ind_{K}^{\QQ_\ell}(\overline{\psi})$
(with $\overline{\psi}$ the reduction of~$\psi$ modulo~$\fP$), is actually isomorphic to
the projectivization of $\rho_{f,\fP}|_{D_\ell}$.
Consequently, $\rhobar_{f,\fP}|_{D_\ell}$ is irreducible as claimed.

Now we consider a prime~$\ell$ satisfying~(\ref{zwei}). In that case we have that
$\widetilde{\rho}_f \cong \Ind_{W_{\ell^2}}^{W_\ell}(\psi)$,
where $\psi$ is a niveau~$2$ character of $\Gal(\Qbar_\ell/\QQ_\ell)$ of order~$n_\ell$.
Hence,
$$ \rhobar_{f,\fP}|_{D_\ell} \cong \Ind_{\QQ_{\ell^2}}^{\QQ_\ell}(\overline{\psi})
\cong \rhobar_{g,\fp}|_{D_\ell}$$
with $\overline{\psi}$ the reduction of~$\psi$ modulo~$\fP$.
Consequently, we find for the Weil-Deligne representation
$(\widetilde{\rho}_g,\widetilde{N}_g)$ of~$g$ that
$\widetilde{\rho}_g \cong \Ind_{W_{\ell^2}}^{W_\ell}(\psi')$ and $\widetilde{N}_g = 0$,
where $\psi'$ is a character of $\Gal(\Qbar_\ell/\QQ_{\ell^2})$
reducing to~$\psi$ modulo~$\fp$.
This means that $\psi' = \psi \alpha$ with some character~$\alpha$ of order
a power of~$p$. As $p$ does not divide $\ell(\ell^2-1)$, it follows that $\alpha$
is unramified and is hence already a character of~$\Gal(\Qbar_\ell/\QQ_\ell)$.
It follows that $\widetilde{\rho}_g \cong \Ind_{W_{\ell^2}}^{W_\ell}(\psi) \otimes \alpha$.
As, however, $f$ and $g$ both have trivial nebentype and the same weigt,
the determinant of $\widetilde{\rho}_g$ is the same as the determinant
of $\widetilde{\rho}_f$ and thus $\alpha^2 = 1$, whence $\alpha$ is the trivial
character. Consequently, $\widetilde{\rho}_f \cong \widetilde{\rho}_g$, proving
that $g$ is tamely dihedral at~$\ell$ of the same order as~$f$.
We have thus actually proved that the Weil-Deligne representations at~$\ell$
of $f$ and $g$ are the same.

Now we consider a prime~$\ell \mid N$ in case~(\ref{drei}).
As $\rho_{f,\fP}|_{D_\ell}$ can be conjugated to lie in the upper triangular matrices,
the same holds for $\rho_{g,\fp}|_{D_\ell}$ by the above discussion on the irreducibility.
Let $(\phi_1,\phi_2)$ and $(\psi_1,\psi_2)$
be the characters on the diagonal of $\rho_{f,\fP}|_{D_\ell}$ and
$\rho_{g,\fp}|_{D_\ell}$, respectively.
As under reduction modulo~$\fP$ (respectively, $\fp$) only $p$-power orders vanish,
it follows from the orders of $\phi_1$ and $\phi_2$ being odd, that this is also
the case for $\psi_1$ and~$\psi_2$.

That $g$ does not have any nontrivial inner twists and no CM again follows from
Theorem~\ref{thm:it}.
\qed

\begin{rem}
At two places in the proof of Theorem~1.1 of~\cite{W} it is implicitly used that the
newforms $f$ and $g$ do not have any nontrivial inner twists, namely
for having $\rhobar_{f,p}(c) \in \PSL_2(\FF_{p^r})$ (in Lemma~3.2)
and for $\rhobar_{g,l}^\proj(\Frob_w) \in \PSL_2(\FF_{l^t})$ (last line but one
of the proof of Theorem~1.1).
This is easily remedied by excluding inner twists
with the help of Proposition~\ref{prop:tamedihedral} and a suitable choice of the starting form~$f$.
If $\ell \not\in \{3,5,7,13\}$, we can just take $f \in S_2(\Gamma_0(\ell))$,
in the other cases we choose $f$ of level a suitable power of~$\ell$
such that $\rho_{f,\iota}(D_\ell)$
can be conjugated into the upper triangular matrices such that the elements on the
diagonal all have odd order.
\end{rem}

\section{Construction of eigenforms without exceptional primes}

The aim of this section is to construct families of Hecke eigenforms without exceptional
primes and without nontrivial inner twists. They also give a uniform way of reproving
(a slightly weaker form of) the main result of~\cite{W}, i.e.\ an application of
modular forms to the inverse Galois problem in the {\em vertical} direction.

\begin{prop}\label{prop:nonexc}
Let $p,q,t,u$ be distinct odd primes and let $N$ be an integer coprime to $pqtu$.
Let $p_1,\dots,p_m$ be the prime divisors of~$6N$.
Let $f \in S_2(Nq^2u^2,\chi)$ be a Hecke eigenform without CM
which is tamely dihedral of order~$p^r>5$ at~$q$ and tamely dihedral of
order~$t^s>5$ at~$u$.
Assume that $q$ and $u$ are completely split in $\QQ(i,\sqrt{p_1},\dots,\sqrt{p_m})$
and that $\legendre q u = 1$ (and, hence, $\legendre u q = 1$ by quadratic reciprocity).

Then $f$ does not have any exceptional primes, i.e.\ for all maximal ideals
$\Lambda$ of~$\cO_{\QQ_f}$, the image of $\rhobar_{f,\Lambda}$ is
$\PSL_2(\FF_\lambda)$ or $\PGL_2(\FF_\lambda)$ in the notation of Theorem~\ref{thm:ribet}.
\end{prop}

\pf
The argument is similar to that used in~\cite{W}, Theorem~1.1, and was inspired by~\cite{KW1}.
Let $\Lambda$ be any maximal ideal of~$\cO_{\QQ_f}$ and
suppose it lies over the rational prime~$\ell$.
Due to the tamely dihedral behaviour, $\rhobar_{f,\Lambda}$ is
irreducible. For, if $\ell \not\in \{p,q\}$, then already $\rhobar_{f,\Lambda}|_{D_q}$
is irreducible. If $\ell \in \{p,q\}$, then $\ell \not\in \{t,u\}$ and already $\rhobar_{f,\Lambda}|_{D_u}$ is irreducible.

We now suppose that the projective image is a dihedral group, i.e.\
it is the induction of a character of a quadratic extension $R/\QQ$, i.e.\
$\rhobar_{f,\Lambda}^\proj \cong \Ind_R^\QQ (\alpha)$ for some character~$\alpha$
of~$\Gal(\Qbar/R)$.
A priori we know from the ramification of $\rhobar_{f,\Lambda}$ that
$R \subseteq \QQ(i,\sqrt{\ell},\sqrt{q},\sqrt{u},\sqrt{p_1},\dots,\sqrt{p_m})$.

Assume first that $\ell \not\in \{p,q\}$. In that case we have
$$ \rhobar^\proj_{f,\Lambda}|_{D_q} \cong \Ind_{\QQ_{q^2}}^{\QQ_q} (\psi)
\cong \Ind_{R_\fQ}^{\QQ_q} (\alpha)$$
for some prime $\cO_R \rhd \fQ \mid q$, where $\psi$ is a niveau~$2$
character of order~$p^r$. From this it follows that $q$ is inert in~$R$.
By assumption, $q$ is split in $\QQ(i,\sqrt{u},\sqrt{p_1},\dots,\sqrt{p_m})$,
whence $R \in \{ \QQ(\sqrt{\ell}), \QQ(\sqrt{-\ell})\}$ and $\ell \nmid 6Nu$.
As $\ell$ is bigger than~$3$ and does not divide the level of~$f$ and the weight is~$2$,
the field $R$ cannot ramify at~$\ell$ either. This was proved by Ribet
(see the proof of Proposition~2.2 in~\cite{Rib97}),
using results of Raynaud implying that the Serre weight is~$2$ and thus
that the projective image of inertia~$I_\ell$ is cyclic of order $\ell+1$ or~$\ell-1$.
We thus obtain a contradiction showing that $\ell \in \{p,q\}$.
In particular, $\ell \not\in \{t,u\}$. Exchanging the
roles $q \leftrightarrow u$, $p \leftrightarrow t$ and $r
\leftrightarrow s$, the very same arguments again lead to a
contradiction. Thus, the projective image of $\rhobar_{f,\Lambda}$
is not a dihedral group.

By the classification of the finite subgroups of $\PGL_2(\Fbar_\ell)$, it remains to exclude
$A_4$, $S_4$ and $A_5$ as projective images. This, however, is clear,
since there is an element of order bigger than~$5$ in the projective image.
\qed

\begin{thm}\label{thm:family}
There exist eigenforms $(f_n)_{n\in \NN}$ of weight~$2$ with trivial nebentype
and without nontrivial inner twists and without CM such that
\begin{enumerate}[(i)]
\item for all~$n$ and all maximal ideals $\Lambda_n \lhd \cO_{f_n}$ the residual Galois
representation $\rhobar_{f_n,\Lambda_n}$ is nonexceptional and
\item for fixed prime~$\ell$, the size of the image of $\rhobar_{f_n,\Lambda_n}$
for $\cO_{\QQ_{f_n}} \rhd \Lambda_n \mid \ell$
is unbounded for running~$n$.
\end{enumerate}
\end{thm}

\pf
Start with some newform $f \in S_2(\Gamma_0(N))$ for squarefree level~$N$.
It does not have any nontrivial inner twists and no CM by Corollary~\ref{cor:noit}.
Let $p_1,\dots,p_m$ be the prime divisors of~$6N$.

Let $B>0$ be any bound.
Let $p$ be any prime bigger than~$B$ provided by Proposition~\ref{prop:tamedihedral}
applied to~$f$ and the set $\{p_1,\dots,p_m\}$, so that we get
$g \in S_2(Nq^2,\chi_\triv)$ which is tamely dihedral at~$q$ of order~$p$
and which does not have any nontrivial inner twists and no CM, for some choice of~$q$.
Now we apply Proposition~\ref{prop:tamedihedral} to~$g$ and the set $\{q,p_1,\dots,p_m\}$,
and get a prime $t>B$ different from~$t$ and some $h \in S_2(Nq^2u^2,\chi_\triv)$
which is tamely dihedral at~$u$ or order~$t$ and which is again without nontrivial inner
twists and without CM, for some choice of~$u$.
By Proposition~\ref{prop:nonexc}, the form~$h$ does not have any
exceptional primes.

We obtain the family $(f_n)_{n \in \NN}$ by increasing the bound~$B$ step by step,
so that elements of bigger and bigger orders appear in the inertia images.
\qed

\begin{rem}\label{rem:vertical}
Theorem~\ref{thm:family} specializes to the following slightly weaker version of
Theorem~1.1 of~\cite{W} concerning the {\em vertical} direction of the inverse Galois problem
for projective linear groups:
\begin{quote}
For every prime~$\ell$, there is an infinite set of natural numbers~$r$
such that $\PSL_2(\FF_{\ell^r})$ or $\PGL_2(\FF_{\ell^r})$ occurs as a Galois group over~$\QQ$.
\end{quote}

The new part here is that the same family $(f_n)_{n \in \NN}$ can be used for all primes~$\ell$.
Theorem~1.1 of~$\cite{W}$ is stronger in the sense that it only concerns $\PSL_2$ and that
it contains a ramification statement.
\end{rem}

\begin{rem}
Note that many choices were made in Theorem~\ref{thm:family} and that one can imagine
many variations in the proof, resulting in many different families.
\end{rem}

\section{Application to the inverse Galois problem}

In this section we prove our main result in the {\em horizontal direction}, i.e.\
Theorem~\ref{thm:main}.
We make use of the only way known to us to impose some condition on the
coefficient field of a newform, namely, by prescribing certain local ramification conditions.
They allow us to have a suitable real cyclotomic field inside the twist invariant
coefficient field.
Such a result is provided by Proposition~\ref{prop:ordp},
when there exists a tamely dihedral prime.
Another such theorem is the following one by Brumer.
In~\cite{D} the first author already observed its usefulness in applications to the
inverse Galois problem.

\begin{thm}[Brumer]\label{thm:brumer}
Let $f \in S_2(N,\chi)$ be a newform without CM.
If $p^{r_p} \mid \mid N$, let $s_p$ be the least integer bigger than or equal to
$$ \frac{r_p}{2} - 1 - \frac{1}{p-1}.$$

Then $\QQ(\zeta_{p^{s_p}}+ \overline{\zeta_{p^{s_p}}})$ is a subfield of~$F_f$.
\end{thm}

\pf
\cite{B}, Theorem 5.5 (b) and the introduction.
\qed
\medskip

We could also give a straightforward argument for $r_p=3$ along the lines of the proof
of Proposition~\ref{prop:ordp}. The existence of a real cyclotomic field inside the
twist invariant coefficient field will allow the application of the following proposition.

\begin{prop}\label{prop:ingredient}
Let $F/\QQ$ be a finite field extension which contains a cyclic field $K/\QQ$ of degree~$n$.
Then the set of primes~$\ell$ such that there is an ideal $\lambda \lhd \cO_F$ dividing~$\ell$
of residue degree~$n$ has a positive density.
\end{prop}

This proposition follows very easily from the following well-known number
theoretic statement.

\begin{prop}\label{prop:galois}
Let $G$ be the Galois group of a normal extension $L/k$ of
number fields, $\ell$ be a prime in $\cO_k$ which is unramified in $\cO_L$,
$\phi$ be any Frobenius automorphism of $\ell$ in~$G$, $H$ be an arbitrary subgroup of $G$,
and $F$ be the subfield of $L$ fixed by $H$. Suppose that the right action of the cyclic
subgroup $\langle \phi \rangle$ of~$G$ partitions the set $H\backslash G$
of right cosets of $H$ into $r$ orbits
with $f_1, \dots , f_r$ cosets respectively.
Then $\ell$ splits into $r$ primes $\lambda_i$ in $\cO_F$, for
which the residual degrees $f(\lambda_i/\ell)$ are given by the numbers $f_i$.
\end{prop}

\pf
\cite{M}, Theorem~33, p.~111.
\qed
\medskip

{\bf Proof of Proposition~\ref{prop:ingredient}.}
Let $L$ be the Galois closure of $F$ over~$\QQ$ with Galois group~$G$.
Let~$H$ be the subgroup of~$G$ such that $F = L^H$.
By assumption, there is a surjection of groups
$G \xrightarrow{\pi} \ZZ/n\ZZ$ such that $H$ is contained in $\ker(\pi)$.

Let $g \in G$ be any element such that $\pi(g)=1 \in \ZZ/n\ZZ$ and
let $r \ge 1$ be the minimum integer such that $g^r \in H$. Then
$ \{ H g^i \;|\; i \in \NN \} \subset H\backslash G$
consists of $r$ elements.
Moreover, as $g^r \in H \subseteq \ker(\pi)$, it follows that
$0 = \pi(g^r)=r \in \ZZ/n\ZZ$, whence $n \mid r$, say, $r = nm$.
Put $\phi=g^m$. Then the minimum $s \ge 1$ such that $\phi^s \in H$ is equal
to~$n$. Consequently, the set
$ \{ H \phi^i \;|\; i \in \NN \} \subset H\backslash G$
consists of~$n$ elements.

By Chebotarev's density theorem, the set of primes~$\ell$ such that
the conjugacy class of the Frobenius elements at~$\ell$ is in the conjugacy class
of $\phi$ in~$G$ has a positive density. By Proposition~\ref{prop:galois},
every such $\ell$ has the property that there is an ideal~$\lambda \lhd \cO_F$ dividing~$\ell$
of residue degree~$n$, as desired.
\qed
\medskip

{\bf Proof of Theorem~\ref{thm:main}}.
We first prove (a) for even~$n$ and (b) for odd~$n$ together. Then we prove (a) for odd~$n$.

We choose a prime~$p$ which is $1 \modulo 2n$.
We also choose an auxiliary prime~$N>13$ different from~$p$.
By Proposition~\ref{prop:tamedihedral} there exists a prime~$q$
(in fact, $q$ can be chosen from a set of primes of positive density)
such that there is a newform $f \in S_2(\Gamma_0(Nq^2))$ without CM
and without nontrivial inner twists which is tamely dihedral of order~$p>1$ at~$q$.
Further, by Proposition~\ref{prop:ordp}, the field $F_f=\QQ_f$ contains
the maximal totally real subfield of the cyclotomic field $\QQ(\zeta_{p})$.

As $p \equiv 1 \modulo 2n$, the number of elements in the group
$\Gal(\QQ(\zeta_{p}+\overline{\zeta_{p}})/\QQ) \cong (\ZZ/p\ZZ)^\times / \{\pm 1\}$
is divisible by~$n$. Hence, there exists a Galois extension $K/\QQ$ with
Galois group $\ZZ/n\ZZ$ such that $K \subseteq F_f=\QQ_f$.

By Proposition~\ref{prop:ingredient}, the set of primes~$\ell$ such that
there is $\lambda \lhd \cO_{F_f}$ dividing~$\ell$ of residue degree~$n$ has a positive density.
For almost all such~$\lambda$,
Ribet's big image theorem (Corollary~\ref{cor:ribet})
now implies that the projective image of the residual Galois representation $\rhobar_{f,\lambda}$
is equal to $\PSL_2(\FF_{\ell^n})$ or $\PGL_2(\FF_{\ell^n})$.

Assuming that $n$ is even, every determinant (an element of $\FF_\ell^\times$) is a square
in $\FF_{\ell^n}$. Consequently, the projective image is equal to $\PSL_2(\FF_{\ell^n})$,
proving (a) for even~$n$.
Assuming, on the other hand, that $n$ is odd, the same reason shows that the projective image is $\PGL_2(\FF_{\ell^n})$, proving~(b).

Now we prove (a) for odd~$n$.
We choose a prime $N \equiv 1 \modulo 2n$.
We also choose an auxiliary prime $q_1 \neq N$ and an auxiliary prime $q_2 \equiv 3 \modulo 4$
different from~$q_1$.
Let $\chi: (\ZZ/q_2\ZZ)^\times \to \{\pm 1\}$ be the unique odd Dirichlet character.
Let $f \in S_3(q_1 q_2 N^3,\chi)$ be a newform (the space is nonempty, as one can
deduce from dimension formulae). The auxiliary prime $q_1$ ensures that
$f$ does not have CM. For, by the definition of the conductor, the image of inertia
$\rho_{f,\Lambda}(I_{q_1})$ contains (after conjugation)
an element of the form $\mat 1b01$ with $b \neq 0$.
Moreover, we know that $[\QQ_f:F_f]$ is a power of~$2$,
as $\chi$ only takes real values (see Section~\ref{sec:it}).

By Theorem~\ref{thm:brumer}, the maximal totally real subfield of $\QQ(\zeta_N)$ is contained
in $F_f$. As above, this implies that there is a field $K \subseteq F_f \subseteq \QQ_f$ such that
$K/\QQ$ is a cyclic extension of order~$n$.
This time we apply Proposition~\ref{prop:ingredient} with the field~$\QQ_f$.
We get that the set of primes~$\ell$ such that
there is $\Lambda \lhd \cO_{\QQ_f}$ dividing~$\ell$ of residue degree~$n$ has a positive density.
Note that the residual degree of $\lambda = \Lambda \cap F_f$ is also equal to~$n$
due to the oddness of~$n$.
For almost all such~$\lambda$,
Ribet's big image theorem (Corollary~\ref{cor:ribet})
yields that the image of the residual Galois representation $\rhobar_{f,\Lambda}$
is equal to $\PSL_2(\FF_{\ell^n})$ or $\PGL_2(\FF_{\ell^n})$.

Note that all determinants are of the form $\pm \FF_\ell^{\times 2}$. Hence, in order
to obtain $\PSL_2(\FF_{\ell^n})$ as projective image, it suffices and we
need to impose $\ell \equiv 1 \modulo 4$.
This is possible. For, if $\QQ(i)$ is disjoint from the Galois closure of $\QQ_f$ over~$\QQ$,
the condition on~$\Lambda$ is independent from $\ell \equiv 1 \modulo 4$. If, however, $\QQ(i)$
is contained in the Galois closure of~$\QQ_f$ over~$\QQ$, then $\QQ_f$ contains~$i$.
As $n$ is odd, any~$\ell$ such that there is $\Lambda \mid \ell$ of odd residue degree
must be $\equiv 1 \modulo 4$.
\qed

\begin{rem}
We point out that in the proof of Theorem~\ref{thm:main} we made many choices.
By varying these choices, the density of the set~$\ell$ such that $\PSL_2(\FF_{\ell^n})$
occurs as a Galois group over~$\QQ$ will certainly increase.
Up to this point, however, we were unable to prove a nontrivial result in this direction.

Moreover, by Weinstein's Theorem~\ref{thm:weinstein},
the auxiliary prime~$N$ in (a) for even~$n$ is not necessary.
\end{rem}

\vspace*{.2cm}

\noindent Luis Dieulefait\\
Departament d'\`Algebra i Geometria\\
Facultat de Matem\`atiques\\
Universitat de Barcelona\\
Gran Via de les Corts Catalanes, 585\\
08007 Barcelona\\
Spain\\
E-mail: {\tt ldieulefait@ub.edu}

\bigskip

\noindent Gabor Wiese\\
Universit\"at Duisburg-Essen\\
Institut f\"ur Experimentelle Mathematik\\
Ellernstra{\ss}e 29\\
45326 Essen\\
Germany\\
E-mail: {\tt gabor.wiese@uni-due.de}\\
Web page: {\tt http://maths.pratum.net/}


\begin{thebibliography}{XXX}

\itemsep=0cm plus 0pt minus 0pt

\bibitem[B]{B} A.~Brumer.
{\it The rank of $J_0(N)$.}
Columbia University Number Theory Seminar (New York, 1992).
Ast\'erisque  No.~228  (1995), 3, 41--68.

\bibitem[DT]{DT} F.~Diamond, R.~Taylor.
{\it Non-optimal levels of mod~$l$ modular representations.}
Invent.\ math.\ {\bf 115} (1994), 435--462.

\bibitem[Di1]{Di1} L.~V.~Dieulefait.
{\it Newforms, Inner Twists, and the Inverse Galois Problem for
Projective Linear Groups}.  J. Th. Nombres Bordeaux {\bf 13} (2001),
395--411.

\bibitem[Di2]{D} L.~V.~Dieulefait.
{\it A control theorem for the images of Galois actions on certain
infinite families of modular forms} in {\it Modular Forms on
Schiermonnikoog}, edited by Gerard van der Geer, Ben Moonen and Bas
Edixhoven, Cambridge University Press, 2008, 79--84.

\bibitem[Di3]{Di3} L.~V.~Dieulefait.
{\it Remarks on Serre's modularity conjecture}. Preprint (2006),
arXiv: math/ 0603439

\bibitem[DV]{DV} L. V. Dieulefait, N. Vila.
{\it Projective linear groups as Galois groups over $\mathbb{Q}$ via
modular representations}. J. Symbolic Comput. {\bf 30} (2000),
799--810.

\bibitem[KLS]{KLS} C.~Khare, M.~Larsen, G.~Savin.
{\it Functoriality and the inverse Galois problem}.
Compos.\ Math.\  {\bf 144}  (2008),  no.\ 3, 541--564.

\bibitem[KLS2]{KLS2} C.~Khare, M.~Larsen, G.~Savin.
{\it Functoriality and the Inverse Galois problem II: groups of type $B_n$ and $G_2$}.
Preprint 2008, arXiv:0807.0861

\bibitem[KW0]{KW0} C.~Khare, J.-P.~Wintenberger.
{\it On Serre's conjecture for 2-dimensional mod~$p$ representations
of $\Gal(\Qbar/\QQ)$}, Annals of Mathematics, {\bf 169} (2009), 229--253.

\bibitem[KW1]{KW1} C.~Khare, J.-P.~Wintenberger.
{\it Serre's modularity conjecture (I)}. Preprint, 2006.

\bibitem[KW2]{KW2} C.~Khare, J.-P.~Wintenberger.
{\it Serre's modularity conjecture (II)}. Preprint, 2007.

\bibitem[Ki]{Ki} M. Kisin.
{\it Modularity of 2-adic Barsotti-Tate representations}. Preprint,
2006.

\bibitem[M]{M} D.~A.~Marcus.
{\it Number Fields}. Universitext. Springer.

\bibitem[Q]{Q} J.~Quer.
{\it Liftings of projective $2$-dimensional Galois representations and embedding problems.}
J.\ Algebra {\bf 171} (1995), no.~2, 541--566.

\bibitem[RV]{RV} A. Reverter, N. Vila.
{\it Some projective linear groups over finite fields as Galois
groups over $\mathbb{Q}$}. Contemp. Math. {\bf 186} (1995), 51--63.

\bibitem[R1]{Ri75}K.~A.~Ribet,
{\it On l-adic representations attached to modular forms}. Invent.
Math. {\bf 28} (1975), 245--275.

\bibitem[R2]{Rib80} K.~A.~Ribet,
{\it Twists of modular forms and endomorphisms of abelian varieties}.
Math.\ Ann.\ {\bf 253} (1980), no.\ 1, 43--62.

\bibitem[R3]{Rib85} K.~A.~Ribet,
{\it On l-adic representations attached to modular forms. II}.
Glasgow Math.\ J.\ {\bf 27} (1985), 185--194.

\bibitem[R4]{Rib97} K.~A.~Ribet,
{\it Images of semistable Galois representations.}
Olga Taussky-Todd: in memoriam. Pacific J.\ Math.\ 1997, Special Issue, 277--297.

\bibitem[T]{T} J.~Tate,
{\it Number Theoretic Background}. Proceedings of Symposia in Pure Mathematics.
Vol.~{\bf 33}, part~2, 3--26.

\bibitem[We]{We} J.~Weinstein.
{\it Hilbert Modular Forms With Prescribed Ramification}. Accepted at Int.\ Math.\ Res.\ Not.\, 2008.

\bibitem[Wi]{W} G.~Wiese.
{\it On projective linear groups over finite fields as Galois groups over the rational numbers}
in {\it Modular Forms on Schiermonnikoog}, edited by
Gerard van der Geer, Ben Moonen and Bas Edixhoven, Cambridge University Press, 2008, 343--350.

\end{thebibliography}
\end{document}